\documentclass{amsart}
\input xy
\xyoption{all}
\usepackage{epsfig}
\usepackage{amsthm}
\usepackage{amssymb}
\usepackage{amsmath}
\usepackage{amscd}

\newtheorem {Theorem}  {Theorem}
\numberwithin{Theorem}{section}
\numberwithin{equation}{section}
\newtheorem {Lemma}[Theorem]  {Lemma}

\theoremstyle{definition}
\newtheorem{Definition}[Theorem]{Definition}
\theoremstyle{remark}
\newtheorem{Remark}[Theorem]{Remark}

\newtheorem {Corollary}[Theorem]{Corollary}

\expandafter\chardef\csname pre amssym.def at\endcsname=\the\catcode`\@ \catcode`\@=11
\def\undefine#1{\let#1\undefined}
\def\newsymbol#1#2#3#4#5{\let\next@\relax
 \ifnum#2=\@ne\let\next@\msafam@\else
 \ifnum#2=\tw@\let\next@\msbfam@\fi\fi
 \mathchardef#1="#3\next@#4#5}
\def\mathhexbox@#1#2#3{\relax
 \ifmmode\mathpalette{}{\m@th\mathchar"#1#2#3}%
 \else\leavevmode\hbox{$\m@th\mathchar"#1#2#3$}\fi}
\def\hexnumber@#1{\ifcase#1 0\or 1\or 2\or 3\or 4\or 5\or 6\or 7\or 8\or
 9\or A\or B\or C\or D\or E\or F\fi}
\font\teneufm=eufm10 \font\seveneufm=eufm7 \font\fiveeufm=eufm5
\newfam\eufmfam
\textfont\eufmfam=\teneufm \scriptfont\eufmfam=\seveneufm
\scriptscriptfont\eufmfam=\fiveeufm

\catcode`\@=\csname pre amssym.def at\endcsname
\newcounter{remark}
\newcommand{\bg}{\begin{equation}}
\newcommand{\ed}{\end{equation}}
\newcommand{\bga}{\begin{eqnarray}}
\newcommand{\eda}{\end{eqnarray}}

\def\cbdu{\hfill{$\Box$}}

\def  \12  {{\frac{1}{2}}}

\begin{document}

\title[The Viscous Camassa-Holm Equations in the Plane]{Decay Asymptotics of the Viscous Camassa-Holm Equations in the Plane}

\author[Clayton Bjorland ]{Clayton Bjorland }

\address{Department of Mathematics, UC Santa Cruz, Santa Cruz, CA 95064,USA}



\thanks{The author was partially supported by NSF Grant OISE-0630623}

\email{cbjorland@math.ucsc.edu} 

\keywords{Camassa-Holm, Navier-Stokes-$\alpha$, vorticity, decay}

\subjclass[2000]{35Q35, 76B03} \

\date{\today}



\bigskip

\begin{abstract}
We consider the vorticity formulation of the 2-D viscous Camassa-Holm equations in the whole space.  We establish global existence for solutions corresponding to initial data in $L^1$ and describe the large time behavior of solutions with sufficiently small and localized initial data.  We calculate the rate at which such solutions approach an ``un-filtered'' Oseen vortex by computing the rate at which the solution of a scaled vorticity problem approaches the solution to a corresponding linearized equation.
\end{abstract}

\maketitle

\section{Introduction}
The viscous Camassa-Holm equations (VCHE) equations are commonly written:
\begin{align}\label{PDE:VCHE}
v_t+u\cdot\nabla v + v\cdot(\nabla u)^T +\nabla p &= \triangle v\\
H^{-1}_\alpha(u) = u-\alpha^2\triangle u &= v\notag\\
\nabla\cdot u &= 0 \notag\\
v(0)&=v_0\notag
\end{align}
Here, $H_\alpha$ is the Helmholtz operator with constant $\alpha$ defined by solving the PDE $u-\alpha^2\triangle u = v$ and will be referred to as the ``filter'' associated with the VCHE.  For a derivation of these equations from variational principles see \cite{MR1627802} or \cite{MR1853633}; for a derivation based on modifying the Navier-Stokes system see \cite{MR1837927}.  The significance of these equations is in the combination of a close relation to the famous Navier-Stokes equations and easily computable bounds on solutions.  In particular, in dimensions 2-4 the VCHE admit smooth global solutions which satisfy the Kelvin Circulation Theorem (circulation is conserved around ``loops'' moving with the filtered flow), the filter is responsible for the smoothing effect on solutions and the non-linear term $v\cdot(\nabla u)^T$ brings the solution into compliance with the Kelvin Circulation Theorem.  These solutions are well suited for numerical and analytic calculations and retain many properties displayed by solutions of the Navier-Stokes equations.  For proofs of global existence and uniqueness in dimensions 2-4 see  \cite{BjorlandSchonbek06}, \cite{MR1837927}, \cite{MR1878243}, \cite{MR2031580}, and \cite{MR1853633}.  Decay of energy and higher norms is considered in \cite{BjorlandSchonbek06}.  Numerical literature is outside the scope of this paper and the reader is referred to \cite{htlans} for a survey of the VCHE role in computational turbulence models and a more complete bibliography.  

The relation between the VCHE and the Navier-Stokes equations is particularly visible in the vorticity form of the equations found by taking the curl ($\nabla\times v = \tilde{v}$) of (\ref{PDE:VCHE}).  In two dimensions the vorticity form is:
\begin{align}\label{PDE:VCHEvorticity}
\tilde{v}_t+u\cdot\nabla \tilde{v} &= \triangle \tilde{v}\\
B(H_\alpha(\tilde{v}))&=u\notag\\
\tilde{v}(0) &=\tilde{v}_0\notag
\end{align}
In the above equation, $B$ represents convolution with the well known Biot-Savart kernel, see \ref{operator:biotsavart}.
The aim of this paper is to further explore the relationship between solutions of the Navier-Stokes equation and the VCHE by describing the way a solution of (\ref{PDE:VCHEvorticity}) approaches the fixed point zero, i.e. computing the first and second order decay asymptotics for solutions with small initial data.    

The asymptotic behavior of the 2-D vorticity equation for the Navier-Stokes equation has been studied in \cite{MR1274542}, \cite{MR1912106} and \cite{MR953819}.  In \cite{MR1912106} the asymptotics are calculated by applying invariant manifold techniques to the semiflow governing the vorticity problem.  Their approach was to scale the vorticity problem into coordinates which are particularly well suited to studying the large time behavior of the Navier-Stokes equation and then apply the Invariant Manifold Theorem in \cite{MR1472350} to construct an invariant manifold in the phase space of the scaled problem and foliate the phase space locally, near the fixed point zero.  The manifolds constructed give insight into the behavior of solutions near the fixed point and, among other results, the authors calculate the asymptotics through the interaction and properties of these manifolds.

The close relation of the Navier-Stokes equations and the VCHE gives hope that a similar program may be carried out for the 2-D VCHE, especially when comparing the vorticity equations.  In fact, such attempts are met with resistance from the filter in the VCHE.  In a functional setting the filter eases problems by smoothing the solution but in a dynamical setting such as this the filter adds complication to the problem.  In particular, the filter does not scale well with the other parts of the equation and the resulting non-linear term has dependence on the scaled time variable not present in the case of the Navier-Stokes equations.  It is not immediately clear how to construct invariant manifolds with such a time dependence or even if such a structure exists.

The Invariant Manifold Theorem in \cite{MR1472350} is based on solving Lyapunov-Perron type equations which are generated through recursive application of the semiflow, the work that follows is based on parts of this approach that work well in the setting of the VCHE.  To construct our decay theorem we linearize the scaled VCHE around a fixed point sufficiently close to zero and determined by the initial data considered, then subtract the linearized equation from the scaled one to get a system which measures distance from the linear solution.  To work around the time dependence of this system we construct an infinite family of systems by steping the original system forward in time a fixed length, these can then be composed to reconstruct the flow of the solution.  Following in spirit \cite{MR1472350}, the semiflow generated by these systems is decomposed into a linear term, a non-linear term, and in the second order case a forcing function which decays quickly for which we can find uniform Lipschitz bounds.  Applying this decomposition recursively we find a discrete Lyapunov-Perron type system which is solvable in a rapidly decaying space.  The existence of a solution to this discrete system implies decay properties of the difference system which in turn allow us to compute the asymptotics we desire.  Although the notion of invariant manifolds is lost we still retain enough structure to complete the asymptotic calculations on the solution.  This procedure, in theory, can be continued to arbitrary orders of asymptotics by considering sufficiently localized data.  We consider only the cases where the governing ODEs are linear.

The work is based in large part on the work in \cite{MR1472350} and \cite{MR1912106}, this is reflected in the notation we use and the statements of many theorems.  Following this introduction we introduce the majority of our notation and a few useful lemmas relating to the VCHE.  Section 4 is dedicated to proving existence and uniqueness of solutions to the vorticity problem in 2-D.  This section is somewhat based on the work in \cite{MR1308857} where similar results were proved for the Navier-Stokes equations but for us the work is simplified thanks to the smoothing properties of the filter in functional settings.  In section 5 we introduce the scaled variables, prove boundedness of solutions with these variables in weighted spaces $L^2(m)$, and discuss other properties of the scaled system.  In this section we also discuss the linearized system and the action of the linear operator $\mathcal{L}$ (the scaled form of $\triangle$) on the weighted spaces.  Sections 6 and 7 hold the theorems and computations involving the first and second order asymptotic respectively.  The theorems in section 7 are very similar to section 6 and many are stated with little or proof when they have nearly identical analogues in section 6.  To end this introduction we state the main conclusions of this paper, proofs are given in sections 4, 6, and 7 respectively.  The spaces $L^2(m)$ are weighted spaces, see (\ref{space:L2m}) below.

\begin{Theorem}\label{theorem:introexistence}
Given initial conditions $\tilde{v}_0\in L^1(\mathbb{R}^2)$, there exists a unique global solution $\tilde{v}\in L^\infty([0,\infty);L^1(\mathbb{R}^2))$ to the PDE (\ref{PDE:VCHEvorticity}).  This solution satisfies the decay bound
\begin{equation}
|\tilde{v}(t)|_q \leq Ct^{-(1-\frac{1}{q})}|\tilde{v}(0)|_1\notag
\end{equation}
\end{Theorem}

\begin{Theorem}\label{theorem:decay1}
For any $\mu\in (0,\frac{1}{2})$, there exists a $r_0>0$ so that for any initial data $\tilde{v}_0\in L^2(2)$ with $\|\tilde{v}_0\|_2\leq r_0$ the solution of (\ref{PDE:VCHEvorticity}) satisfies
\begin{equation}
|\tilde{v}(\cdot,t)-a(\Omega(\cdot,t)-\alpha^2\triangle\Omega(\cdot,t))|_p\leq C(1+t)^{-1-\mu+\frac{1}{p}}\notag
\end{equation}
where $a = \int_{\mathbb{R}^2} \tilde{v}_0\,dx$ and
\begin{equation}
\Omega(x,t)= \frac{1}{4\pi(1+t)}e^{\frac{-|x|^2}{4(1+t)}}\notag
\end{equation}
\end{Theorem}

\begin{Theorem}\label{theorem:decay2}
For any $\mu\in (\frac{1}{2},1)$, there exists a $r_0>0$ so that for any initial data $\tilde{v}_0\in L^2(3)$ with $\|\tilde{v}_0\|_3\leq r_0$ the solution of (\ref{PDE:VCHEvorticity}) satisfies
\begin{equation}
|\tilde{v}(\cdot,t)-a(\Omega(\cdot,t)-\alpha^2\triangle\Omega(\cdot,t))-\sum_{i=1,2}b_i(\partial_i\Omega(x,t)-\alpha^2\triangle\partial_i\Omega(x,t))|_p\leq C(1+t)^{-1-\mu+\frac{1}{p}}\notag
\end{equation}
where $b_i=\int\xi_i\tilde{w}\,d\xi$ and $a$, $\Omega$ are as in the previous theorem.
\end{Theorem}

\section{Notation}
Throughout this paper we use $\mathbb{N}$ to refer to the natural numbers including $0$.  Standard Lebesgue spaces will be denoted $L^p(\mathbb{R}^2)$ or ($L^p$ for short) with the norm $|\cdot|_p=\int \cdot\,dx$.   

To denote the curl of a vector field, we use the tilde.  For example, the curl of a vector field $v$ is given by $\nabla\times v=\tilde{v}$.  For a divergence free vector field, the curl can be undone through convolution with the well known Biot-Savart kernel  $x^\perp/(2\pi|x|^2)$, $x^\perp = (-x_2,x_1)^T$.  We denote this convolution as 
\begin{equation}\label{operator:biotsavart}
B(\tilde{w})=\frac{1}{2\pi}\int_{\mathbb{R}^2}\frac{(x-y)^\perp}{|x-y|^2}\tilde{w}(y)\,dy
\end{equation}
so that $B(\tilde{v})=v$ for divergence free vector fields.

We define $b(\xi)=(1+|\xi|^2)^{1/2}$ and make frequent use of the weighted Hilbert spaces $L^2(m)$ defined by
\begin{align}
L^2(m)&=\{f\in L^2(\mathbb{R}^2)|\ \|f\|_m<\infty\ \&\ \nabla\cdot f =0\}\label{space:L2m}\\
\|f\|^2_m&=\int_{\mathbb{R}^2}b^{2m}(\xi)|f(\xi)|^2\,d\xi\notag
\end{align}

\section{Preliminaries}
This section contains bounds that will be useful later in the paper, particularly section 3.  To begin we recall a lemma from \cite{MR1912106} concerning the Biot-Savart operator $B$ defined by (\ref{operator:biotsavart}), and the curl of divergence free vector fields.
\begin{Lemma}\label{lemma:biotsavart}
Let 
\begin{equation}
v = B(\tilde{v}) = \frac{1}{2\pi}\int_{\mathbb{R}^2} \frac{(x-y)^\perp}{|x-y|^2}\tilde{v}(y)\,dy\notag
\end{equation}

(i)  If $1<p<2<q<\infty$, $\frac{1}{q}=\frac{1}{p}-\frac{1}{2}$, and $\tilde{v}\in L^p(\mathbb{R}^2$, then there exists a $C>0$ such that $|v|_q\leq C|\tilde{v}|_p$.

(ii)  If $1\leq p<2<q \leq\infty$, $\frac{1}{2} =\frac{\alpha}{p} + \frac{(1-\alpha)}{q}$ where $\alpha\in(0,1)$, and $\tilde{v}\in L^p\cap L^q(\mathbb{R}^2$, then there exists a $C>0$ such that $|v|_\infty\leq C|\tilde{v}|^\alpha_p|\tilde{v}|^{1-\alpha}_q$.

(iii)  If $1<p<\infty$ and $\tilde{v}\in L^p(\mathbb{R}^2)$, then there exists a $C>0$ such that $|\nabla v|_p\leq C|\tilde{v}|_p$
\end{Lemma}
\begin{proof}
This is Lemma 2.1 in \cite{MR1912106}.
\end{proof}

The following bounds for the heat kernel will also be useful in proving existence of solutions in section 3.
\begin{Lemma}\label{lemma:heat}
Let $\Phi$ be the fundamental solution to the heat equation.  Then
\begin{align}
|\Phi(t)|_p&\leq\frac{C(n,p)}{ t^{(1-\frac{1}{p})\frac{n}{2}}} \notag\\
|\nabla\Phi(t)|_p &\leq \frac{C(n,p)}{t^{\frac{1}{2}+(1-\frac{1}{p})\frac{n}{2}}}\notag
\end{align}
\end{Lemma}
\begin{proof}
By direct calculation
\begin{equation}
|\Phi(t)|_p^p=\frac{1}{(4\pi t)^{\frac{pn}{2}}}\int_{\mathbb{R}^n}e^{\frac{-p|x|^2}{4t}}\,dx=\frac{1}{p^{\frac{n}{2}}(4\pi t)^{\frac{(p-1)n}{2}}}\notag
\end{equation}
This proves the first bound.  For the second, start by differentiating the heat kernel then take the $L^p$ norm
\begin{equation}
|\nabla\Phi(t)|_p^p \leq \frac{C(n,p)}{t^{p(1+\frac{n}{2})}}\int_{\mathbb{R}^n} |x|^pe^{\frac{-p|x|^2}{4t}}\, dx\notag
\end{equation}
Now a change of variables and integration proves the second bound.
\end{proof}

We also recall some facts about the Helmholtz operator $H$.

\begin{Lemma}
Given, $v\in C^\infty_0(\mathbb{R}^n)$ let $u$ be the smooth solution to the Helmholtz equation
\begin{equation}
H_\alpha^{-1}(u) = u-\alpha^2\triangle u = v\notag
\end{equation}
Then $u$ satisfies the following bounds
\begin{align}
|u|_p&\leq |v|_p \ for \ p\in [1,\infty]\notag\\
|u|_p &\leq \frac{C(n,p,r)}{\alpha^{1+\gamma}}|v|_r \ for \ \gamma=\frac{n}{2}(\frac{1}{r}-\frac{1}{p})<1\notag\\
|\nabla u|_p &\leq \frac{C(n,p,r)}{\alpha^{\frac{3}{2}+\gamma}} |v|_r \ for \ \gamma =\frac{n}{2}(\frac{1}{r}-\frac{1}{p})<\frac{1}{2}\notag
\end{align}
\end{Lemma}
\begin{proof}
See \cite{BjorlandSchonbek06}.
\end{proof}

\begin{Lemma}\label{lemma:helmholtz}
Let $v\in L^r(\mathbb{R}^2)$, $r\in [1,\infty)$.  Then, there exists a solution $u\in W^{1,r}$ to the Helmholtz equation
\begin{equation}
H_\alpha^{-1}(u) = u-\alpha^2\triangle u = v\notag
\end{equation}
which satisfies the bounds
\begin{align}
|u|_p&\leq |v|_p \ for \ p\in [1,\infty]\notag\\
|u|_p &\leq \frac{C(n,p,r)}{\alpha^{1+\gamma}}|v|_r \ for \ \gamma=\frac{n}{2}(\frac{1}{r}-\frac{1}{p})<1\notag\\
|\nabla u|_p &\leq \frac{C(n,p,r)}{\alpha^{\frac{3}{2}+\gamma}} |v|_r \ for \ \gamma =\frac{n}{2}(\frac{1}{r}-\frac{1}{p})<\frac{1}{2}\notag
\end{align}
\end{Lemma}
\begin{proof}
Let $v_n\in C_0^\infty$ be a sequence such that $v_n\rightarrow v$ in $L^r(\mathbb{R}^2)$ and $\|v_n\|_r\leq \|v\|_r$.  By the previous lemma, there exists a sequence of solutions to the Helmholtz equation $u_n$, each term satisfying the desired bounds.  Applying the Banach-Alaoglu theorem to a possible subsequence of the $u_n$ gives the existence of a function $u$ with the desired properties.
\end{proof}

The last Lemma in this section concerns the inclusion $L^2(m)\hookrightarrow L^1$ which will be useful in proving many estimates later in this paper.
\begin{Lemma}
If $f\in L^2(m)$ for some $m>1$ then $f\in L^q$ for $1\leq q \leq 2$.
\end{Lemma}
\begin{proof}
The bound $|f|_2\leq \|f\|_m$ is immediate from the definition of $L^2(m)$.  Since $(1+|\xi|^2)^{\frac{m}{2}}$ is integrable in $\mathbb{R}^2$ if $m>1$ we have $|f|_1\leq C\|f\|_m$.  Interpolation finishes the proof.
\end{proof}

\section{Vorticity Problem for the VCHE}
This section contains proofs for existence and uniqueness of solutions to the vorticity form of the VCHE (\ref{PDE:VCHEvorticity}) with data in $\tilde{v}\in L^1(\mathbb{R}^2)$.  Decay of these solutions is provided by the optimal smoothing results in \cite{MR1381974}.  

Instead of working directly with (\ref{PDE:VCHEvorticity}) we prove existence and uniqueness for the mild form of the PDE (\ref{PDE:VCHEvorticity}) by solving the integral equation
\begin{align}\label{mildPDE:VCHEvorticity}
\tilde{v}(t)&=e^{\triangle t}\tilde{v}_0 - \int_0^t\nabla\Phi(t-s)\ast [B(H_\alpha(\tilde{v}))\otimes\tilde{v}](s)\, ds\\
\tilde{v}(0)&=\tilde{v}_0\notag
\end{align}

The corresponding result for the Navier-Stokes equation was proved in \cite{MR1308857} using a fixed point argument in a space smaller then $L^1$, then extending the solution operator to $L^1$, see also \cite{MR1308858} and \cite{MR1270113}.  This approach was necessary because of the difficulty in writing the Navier-Stokes equations so a fixed point argument can be applied in $L^1$.  In the case of the VCHE, the filter provides enough leverage to apply a fixed point argument directly.  We start by proving a bound on the bilinear term then follow with existence of a global solution.

\begin{Lemma}\label{lemma:bilinearbound}
The bilinear form 
\begin{equation}
\mathcal{B}(\tilde{v},\tilde{w}):L^\infty([0,T],L^1(\mathbb{R}^2))\times L^\infty([0,T],L^1(\mathbb{R}^2))\rightarrow L^\infty([0,T],L^1(\mathbb{R}^2))\notag
\end{equation}
defined by
\begin{equation}
\mathcal{B}(\tilde{v},\tilde{w})=\int_0^t\nabla\Phi(t-s)\ast [B(H_\alpha(\tilde{v}))\otimes\tilde{w}](s)\, ds\notag
\end{equation}
satisfies the bound
\begin{equation}
\sup_{t\in[0,T]}|\mathcal{B}(\tilde{v},\tilde{w})(t)|_1\leq C(T)\left(\sup_{t\in[0,T]}|v(t)|_1\right)\left(\sup_{t\in[0,T]}|w(t)|_1\right)\notag
\end{equation}
where $C(T)\rightarrow 0$ as $T\rightarrow 0$.
\end{Lemma}
\begin{proof}
Apply first Young's inequality then H\"{o}lder's inequality to the bilinear form with $1<p<2<q<\infty$ and $\frac{1}{p}+\frac{1}{q}=1$:
\begin{align}
|B(\tilde{v},\tilde{w})(t)|_1&\leq \int_0^t|\nabla\Phi(t-s)|_1|B(H_\alpha(\tilde{v}))(s)|_\infty|\tilde{w}(s)|_1\, ds\notag
\end{align}
Lemmas \ref{lemma:biotsavart} and \ref{lemma:helmholtz} give the bound $|B(H_\alpha(\tilde{v}))(s)|_\infty\leq C |\tilde{v}(s)|_1$ so
\begin{equation}
|B(\tilde{v},\tilde{w})(t)|_1\leq C\int_0^t|\nabla\Phi(t-s)|_1	\, ds\left(\sup_{s\in[0,t]}|\tilde{v}(t)|_1\right)\left(\sup_{s\in[0,t]}|\tilde{w}(t)|_1\right)\notag
\end{equation}
By Lemma \ref{lemma:heat}, $|\nabla\Phi(t-s)|_1$ is integrable over finite intervals.  Taking the supremum over $t\in[0,T]$ yields
\begin{equation}
\sup_{t\in[0,T]}|B(\tilde{v},\tilde{w})(t)|_1 \leq CT^{\frac{1}{2}}\left(\sup_{t\in[0,T]}|\tilde{v}(t)|_1\right)\left(\sup_{t\in[0,T]}|\tilde{w}(t)|_1\right) \notag
\end{equation}
This estimate concludes the proof.
\end{proof}

\begin{Theorem}\label{PDE:VCHEvorticityexistence}
Given initial data $\tilde{v}_0\in L^1(\mathbb{R}^2)$ there exists a unique global solution $\tilde{v}\in L^\infty([0,\infty);L^1(\mathbb{R}^2))$ to the integral equation (\ref{mildPDE:VCHEvorticity}).  For any $p\in [1,\infty]$, this solution satisfies the decay bound
\begin{equation}\label{opsmoothingbound}
|\tilde{v}(t)|_p \leq Ct^{-(1-1/p)}|\tilde{v}(0)|_1
\end{equation}
\end{Theorem}
\begin{proof}
Lemma \ref{lemma:bilinearbound} with a standard fixed point argument gives the existence of a mild solution in some possibly small time interval.  The length of the time interval depends on the $L^1$ norm of the initial data.  After applying the optimal smoothing results in \cite{MR1381974} (see also \cite{MR916761}) the $L^1$ norm of the solution does not increase beyond the $L^1$ norm of the data which implies the existence of a global solution.  Indeed, the vorticity equation for the VCHE is a viscously damped conservation law so after applying Theorem 1 in \cite{MR1381974} we establish (\ref{opsmoothingbound}).  In particular, $|v(t)|_1\leq C|v(0)|_1$, which establishes global existence.

It remains to establish uniqueness.  Let $\tilde{v}$ and $\tilde{w}$ be two solutions of (\ref{mildPDE:VCHEvorticity}) corresponding to the same initial data $\tilde{v}_0\in L^1(\mathbb{R}^2)$.  After adding and subtracting cross terms we see
\begin{align}
(\tilde{v}-\tilde{w})(t)= -&\int_0^t\nabla\Phi(t-s)\ast[B(H(\tilde{v},\alpha))\otimes(\tilde{v}-\tilde{w})](s)\,ds \notag\\
&+\int_0^t\nabla\Phi(t-s)\ast[B(H((\tilde{v}-\tilde{w}),\alpha))\otimes\tilde{w}](s)\,ds \notag
\end{align}
Similar to the proof of Lemma \ref{lemma:bilinearbound}, apply Lemmas \ref{lemma:biotsavart} and \ref{lemma:helmholtz} with Young's inequality, for $1<p<2$,
\begin{equation}
|\tilde{v}(t)-\tilde{w}(t)|_1 \leq C \int_0^t|\nabla\Phi(t-s)|_1(|\tilde{v}(s)|_1+|\tilde{w}(s)|_1)|\tilde{v}(s)-\tilde{w}(s)|_1\,ds\notag
\end{equation}
Apply (\ref{opsmoothingbound}) to obtain
\begin{equation}
|\tilde{v}(t)-\tilde{w}(t)|_1 \leq C |\tilde{v}_0|_1\int_0^t|\nabla\Phi(t-s)|_1|\tilde{v}(s)-\tilde{w}(s)|_1\,ds\notag
\end{equation}
From here, the Gronwall inequality with Lemma \ref{lemma:heat} is used to establish uniqueness and conclude the proof.
\end{proof}

This establishes Theorem \ref{theorem:introexistence}.

\section{The Scaled Equations}
In this section we introduce scaled variables and rewrite the vorticity equations for the VCHE in these variables.  An existence theorems for the the scaled VCHE and related filter equations in the spaces $L^2(m)$ is the provided.  Finally we discuss the action of the linear operator $\mathcal{L}$ on the weighted spaces $L^2(m)$.  

The scaled variables are defined
\begin{align}
\xi&=\frac{x}{\sqrt{1+t}} &\tau&=\ln(1+t)\notag\\
v(x,t)&=\frac{1}{\sqrt{1+t}}w(\xi,\tau) &u(x,t)&=\frac{1}{\sqrt{1+t}}\omega(\xi,\tau)\notag\\
\tilde{v}(x,t)&=\frac{1}{1+t}\tilde{w}(\xi,\tau) &\tilde{u}(x,t)&=\frac{1}{1+t}\tilde{\omega}(\xi,\tau)\notag
\end{align}
It has been show in \cite{MR1274542}, \cite{MR1912106} and \cite{MR953819} that these variables are very useful when studying the large time behavior of the Navier-Stokes equation.  
The VCHE (\ref{PDE:VCHEvorticity}) becomes
\begin{align}
\tilde{w}_\tau&=\mathcal{L}\tilde{w} -\omega\cdot \nabla_\xi \tilde{w} &\tilde{w}(0)&=\tilde{w}_0\label{PDE:VCHEscaledvorticity}\\
\omega &=B(\mathcal{H}_{\alpha,\tau}(\tilde{w}))\label{PDE:VCHEscaledvorticityfilter}
\end{align}
where $B$ is again convolution with the Biot-Savart kernel as in (\ref{operator:biotsavart}), $\mathcal{L}$ is the linear operator
$\mathcal{L}=\triangle_\xi +\frac{1}{2}\xi\cdot\nabla_\xi+ I$, and  $\mathcal{H}_{\alpha,\tau}$ is the operator defined by solving the scaled Helmholtz equations
$\tilde{\omega} -\alpha^2e^{-\tau}\triangle_\xi \tilde{\omega} = \tilde{w}$.

Our first goal of the section will be to show how the filter $\mathcal{H}_{\alpha,\tau}$ acts on the weighted spaces $L^2(m)$, in particular we will show that the Helmholtz equation has a unique solution on these spaces.

\begin{Lemma}\label{lemma:Helmholtzscaledbound}
If $w\in L^2(m)$ and $\omega\in L^2(m)$ are related by the scaled Helmholtz equation $w = \omega-\alpha^2e^{-\tau}\triangle_\xi \omega$ then $\|\omega(\tau)\|_m^2\leq C\|w(\tau)\|_m^2$.  In the case $4m^2\alpha^2<1$, $C=1$, otherwise $C=2(1+2^{m-1}(4m^2\alpha^2)^m$).
This lemma shows how the operator $\mathcal{H}_{\alpha,\tau}:L^2(m)\rightarrow L^2(m)$ is bounded.
\end{Lemma}
\begin{proof}
The case of $m=0$ follows from the linear theory, the bound is
\begin{equation}
\|w\|^2_0=\|\omega\|^2_0+\alpha^2e^{-\tau}\|\nabla\omega\|_0^2+\alpha^4e^{-2\tau}\|\triangle\omega\|_0^2\notag
\end{equation}
If $m>0$, square the PDE and multiply by $b^{2m}$, after integration we have
\begin{equation}
\|w\|_m^2 = \|\omega\|_m^2+\alpha^4e^{-2\tau}\|\triangle\omega\|_m^2
-2\alpha^2e^{-\tau}\int_{\mathbb{R}^2}b^{2m}\omega\triangle\omega\, d\xi\notag
\end{equation}
Using integration by parts
\begin{equation}
\int_{\mathbb{R}^2}(1+|\xi|^2)^m\omega\triangle\omega\, d\xi=-\|\nabla\omega\|_m^2
+2m\int_{\mathbb{R}^2}(b^{2m-2}+(m-1)b^{2m-4}|\xi|^2)\omega^2\,d\xi\notag
\end{equation}
leaves
\begin{align}
\|w\|_m^2 = \|\omega\|_m^2&+\alpha^4e^{-2\tau}\|\triangle\omega\|_m^2
+2\alpha^2e^{-\tau}\|\nabla\omega\|_m^2\notag\\
&-4m\alpha^2e^{-\tau}\|\omega\|_{m-1}^2-4m(m-1)\alpha^2e^{-\tau}\||\xi|\omega\|_{m-2}^2\notag
\end{align}
The bound $-e^{-\tau}b^{2m-4}(b^2+(m-1)|\xi|^2)\geq -m2b^{2m}$ shows
\begin{equation}\label{weighted:energy}
\|w\|_m^2\geq (1-4m^2\alpha^2)\|\omega\|_m^2+\alpha^4e^{-2\tau}\|\triangle\omega\|_m^2
+2\alpha^2e^{-\tau}\|\nabla\omega\|_m^2\notag
\end{equation}
and proves the result if $4m^2\alpha^2<1$.  

If $4m^2\alpha^2\geq 1$, set $\beta^2=8m^2\alpha^2-1$ so that for $|\xi|\geq\beta$ we have
$b^{2m}\geq 8m^2 \alpha^2 b^{2m-2}$.  If $B(\beta)$ is the ball with radius $\beta$,
\begin{equation}
-4m^2\alpha^2\int_{\mathbb{R}^2}b^{2m-2}\omega^2\,d\xi\geq -4m^2\alpha^2(8m^2\alpha^2)^{m-1}\int_{B(\beta)}\omega^2\,d\xi
-\frac{1}{2}\int_{B^C(\beta)}b^{2m}\omega^2\,d\xi\notag
\end{equation}
Applying the case $m=0$ and the bound $\|\tilde{w}\|_2^2\leq \|\tilde{w}\|_m^2$ allows 
\begin{equation}
-\int_{B(\beta)}\omega^2\,d\xi\geq-\int_{\mathbb{R}^2}\omega^2\,d\xi\geq-\int_{\mathbb{R}^2}w^2\,d\xi\geq -\|w\|^2_m\notag
\end{equation}
Also,
\begin{equation}
-\frac{1}{2}\int_{B^C(\beta)}b^{2m}\omega^2\,d\xi\geq -\frac{1}{2}\|\omega\|_m^2\notag
\end{equation}
Considering all of this with leaves,  in the case $4m^2\alpha^2\geq 1$,
\begin{equation}
C\|w\|_m^2 \geq \frac{1}{2}\|\omega\|_m^2+\alpha^4e^{-2\tau}\|\triangle\omega\|_m^2
+2\alpha^2e^{-\tau}\|\nabla\omega\|_m^2\label{bound:Helmholtzscaledbound}
\end{equation}
where $C=1+2^{m-1}(4m^2\alpha^2)^m$.
\end{proof}
The above proof is some indication that the filter is not well suited for the weighted spaces.  There is still a smoothing effect, but as can be seen from (\ref{bound:Helmholtzscaledbound}) this effect decreases as $\tau$ becomes large.  

\begin{Theorem}\label{lemma:Helmholtzscaledexistence}
Given $w\in L^2(m)$ there exists a unique solution $\omega\in L^2(m)$ to the scaled Helmholtz equations $w = \omega-\alpha^2e^{-\tau}\triangle_\xi \omega$.  This proves the operator $\mathcal{H}_{\alpha,\tau}:L^2(m)\rightarrow L^2(m)$ is well defined.
\end{Theorem}
\begin{proof}
The rough estimate $b^{2m}>1$ implies $L^2(m)\subset L^2(0)$, so if $w\in L^2(m)$ it is well known that there is a unique $\omega\in L^2(0)$ solving the equation.  Lemma \ref{lemma:Helmholtzscaledbound} shows that $\omega\in L^2(m)$.
\end{proof}

We now turn our attention to the scaled VCHE.  Using the strongly continuous semigroup $e^{\tau\mathcal{L}}$ generated by $\mathcal{L}$ in $L^2(m)$ we write the mild form of the scaled vorticity problem as
\begin{align}
\tilde{w}(\tau) &=e^{\tau\mathcal{L}}\tilde{w}_0 -\int_0^\tau e^{-\frac{1}{2}(\tau-s)}\nabla\cdot e^{(\tau-s)\mathcal{L}}(\omega(s)\tilde{w}(s))\,ds\label{mildPDE:VCHEscaledvorticity}\\
\tilde{\omega}&= \mathcal{H}_{\alpha,\tau}(\tilde{w}) \ \ \ \ \ \ \ \tilde{w}(0)=\tilde{w}_0 \notag
\end{align}

The following lemma provides an estimate on the bilinear form which will be used to prove existence of a local solution to (\ref{mildPDE:VCHEscaledvorticity}).
\begin{Lemma}\label{lemma:mbilinear}
Given $\tilde{w}_1,\tilde{w}_2\in C^0([0,T];L^2(m))$, define
\begin{equation}
R(\tilde{w}_1,\tilde{w}_2)(\tau)=\int_0^\tau e^{-\frac{1}{2}(\tau-s)}\nabla\cdot e^{(\tau-s)\mathcal{L}}(\omega_1(s)\tilde{w}_2(s))\,ds\notag
\end{equation}
where $\omega_1=w_1-\alpha^2e^{-\tau}\triangle w_1$ and $w_1$ is obtained from $\tilde{w}_1$ via the Biot-Savart law.  Then $R\in C^0([0,T],L^2(m))$, and there exists $C_0=C_0(m,T)>0$ such that
\begin{equation}
\sup_{0\leq\tau\leq T}\|R(\tilde{w}_1,\tilde{w}_2)(\tau)\|_m\leq C_0\left(\sup_{0\leq\tau\leq T}\|\tilde{w}_1(\tau)\|_m\right)\left(\sup_{0\leq\tau\leq T}\|\tilde{w}_2(\tau)\|_m\right)\notag
\end{equation}
Moreover, the constant $C_0$ becomes arbitrarily small as $T$ tends to zero.
\end{Lemma}
\begin{proof}
We rely on an estimate of the semigroup $e^{\tau\mathcal{L}}$ proved in the appendix of \cite{MR1912106}, if $r\in(\frac{2}{m+1},2)$ then
\begin{equation}
|b^m\nabla \cdot e^{\tau\mathcal{L}} u|_2\leq \frac{C}{a(\tau)^{(\frac{1}{r}-\frac{1}{2})+\frac{1}{2}}}|b^m u|_r\notag
\end{equation}
where $a(\tau)=1-e^{-\tau}$.  This allows
\begin{equation}
|b^m R(\tilde{w}_1,\tilde{w}_2)(\tau)|_2\leq C\int^\tau_0\frac{1}{a(\tau-s)^{\frac{1}{r}}}|b^m \omega_1(s)\tilde{w}_2(s)|_r \,ds\notag
\end{equation}
H\"{o}lder's inequality, Lemma \ref{lemma:biotsavart}, and the inclusion $L^2(m)\hookrightarrow L^r(\mathbb{R}^2)$ provides the bound
$|b^m \omega_1(s)\tilde{w}_2(s)|_r\leq C\|\tilde{w}_2\|_m\|\tilde{\omega}_1\|_m$.
Apply Lemma \ref{lemma:Helmholtzscaledbound} and note that $a(\tau-s)^{-1/r}$ is integrable from $0$ to $\tau$ to finish the proof.
\end{proof}

In addition to providing existence and uniqueness, the following theorem shows how we can control the $L^2(m)$ norm of a solution to the integral equation (\ref{mildPDE:VCHEscaledvorticity}) by controlling the $L^2(m)$ norm of the initial data.  Just as Lemma \ref{lemma:mbilinear} has an analogue in \cite{MR1912106}, our existence proof is modeled after the proof of Theorem 3.2 in \cite{MR1912106}.

\begin{Theorem}\label{theorem:VCHEscaledvorticityexistence}
Given $\tilde{w}_0\in L^2(m)$ for some $m>1$, there exists a global solution $\tilde{w}\in C^0([0,\infty),L^2(m))$ to the integral equation (\ref{mildPDE:VCHEscaledvorticity}) with $\tilde{w}(0)=\tilde{w}_0$.  Moreover, there exists a constant $C_1=C_1(\|\tilde{w}_0\|_m)$ such that
\begin{equation}\label{theorem:VCHEscaledvorticitybound}
\|\tilde{w}(\tau)\|_m\leq C_1
\end{equation}
and $C_1 \rightarrow 0$ as $\|\tilde{w}_0\|_m\rightarrow 0$.
\end{Theorem}
\begin{proof}
The previous lemma and a fixed point argument gives unique local in time existence of a unique solution.  Moreover, there exists a $T>0$ such that
\begin{equation}\label{globalexistencebound0}
\sup_{0\leq\tau\leq T}\|\tilde{w}(\tau)\|_m\leq 2\|\tilde{w}_0\|_m
\end{equation}
By scaling (\ref{opsmoothingbound}) and using the fact that $L^2(m)\hookrightarrow L^1$ for $m>1$ we see that this solution satisfies, for all $p\in [1,\infty]$,
\begin{equation}\label{globalexistencebound1}
|\tilde{w}(\tau)|_p \leq \frac{C_p\|\tilde{w}(0)\|_m}{a(\tau)^{1-\frac{1}{p}}}
\end{equation}

We will now establish (\ref{theorem:VCHEscaledvorticitybound}) which will imply global in time existence.  Multiplying (\ref{PDE:VCHEscaledvorticity}) by $b^{2m}\tilde{w}$ and integrating we find
\begin{equation}
\frac{1}{2}\frac{d}{d\tau}\int_{\mathbb{R}^2} b^{2m} \tilde{w}^2 \,d\xi = \int_{\mathbb{R}^2} b^{2m}\left( \tilde{w}\triangle\tilde{w} +\frac{\tilde{w}}{2}(\xi\cdot\nabla)\tilde{w}
+\tilde{w}^2-\tilde{w}(\omega\cdot\nabla)\tilde{w}\right) \,d\xi\notag
\end{equation}
Integration by parts and the bound $|\xi|\leq b(\xi)$ give the following estimates:
\begin{align}
\int_{\mathbb{R}^2} b^{2m}\tilde{w}\triangle\tilde{w}\,d\xi &\leq -\int_{\mathbb{R}^2} b^{2m}\nabla\tilde{w}^2\,d\xi 
+ 2m^2\int_{\mathbb{R}^2} b^{2m-2}\tilde{w}^2\,d\xi\notag\\ 
\frac{1}{2}\int_{\mathbb{R}^2} b^{2m}\tilde{w}(\xi\cdot\nabla \tilde{w}) 
&= -\frac{1}{2}\int_{\mathbb{R}^2} b^{2m} \tilde{w}^2\,d\xi-\frac{1}{2}m\int_{\mathbb{R}^2} b^{2m-2}|\xi|^2\tilde{w}^2\,d\xi\notag\\
-\int_{\mathbb{R}^2} b^{2m} \tilde{w}(\omega\cdot\nabla)\tilde{w}\,d\xi 
&= m\int_{\mathbb{R}^2} b^{2m-2}(\xi\cdot\omega)\tilde{w}^2\,d\xi \notag
\end{align}
Furthermore, given $\epsilon>0$ there exists a $C_\epsilon >0$ such that $b^{2m-2}\leq \epsilon b^{2m} + C_\epsilon$ so we can bound
\begin{align}
2m^2\int_{\mathbb{R}^2} b^{2m-2}\tilde{w}^2\,d\xi\leq \epsilon \int_{\mathbb{R}^2} b^{2m}\tilde{w}^2\,d\xi 
+ C_\epsilon \int_{\mathbb{R}^2} \tilde{w}^2\,d\xi \notag
\end{align}
Similarly,
\begin{align}
m\int_{\mathbb{R}^2} b^{2m-2}(\xi\cdot\omega)\tilde{w}^2\,d\xi \leq \epsilon \int_{\mathbb{R}^2} b^{2m}\tilde{w}^2\,d\xi 
+ C_\epsilon |\omega|^{2m}_\infty \int_{\mathbb{R}^2} \tilde{w}^2\,d\xi \notag
\end{align}
Putting these bounds together yields
\begin{align}
\frac{1}{2}\frac{d}{d\tau}\int_{\mathbb{R}^2} b^{2m} \tilde{w}^2 \,d\xi \leq -(\frac{1}{2}-2\epsilon)\int_{\mathbb{R}^2} b^{2m} \tilde{w}^2\,d\xi
+C_\epsilon (1+|\omega|_\infty^{2m})\int_{\mathbb{R}^2} \tilde{w}^2\,d\xi\notag
\end{align}
Write $\delta=\frac{1}{2}-4\epsilon$, then
\begin{align}
\frac{d}{d\tau}\left(e^{\delta\tau}\int_{\mathbb{R}^2}b^{2m}\tilde{w}^2\,d\xi\right) \leq C_\epsilon e^{\delta\tau}(1+|\omega|_\infty^{2m})\int_{\mathbb{R}^2} \tilde{w}^2\,d\xi\notag
\end{align}
This inequality implies (\ref{theorem:VCHEscaledvorticitybound}) provided 
\begin{equation}\label{globalexistencebound2}
\sup_{0\leq \tau\leq \infty}\left((1+|\omega|_\infty^{2m})\int_{\mathbb{R}^2} \tilde{w}^2\,d\xi\right)\leq C(\|\tilde{w}_0\|_m)
\end{equation}
where $C(\|\tilde{w}_0\|_m)\rightarrow 0$ as $\|\tilde{w}_m\|_m\rightarrow 0$.  We now establish this estimate.  

When $T$ is such that (\ref{globalexistencebound0}) holds, estimate (\ref{globalexistencebound1}) with $p=\infty$ and then $p=2$ implies (\ref{globalexistencebound2}) when the supremum is taken over $\tau\in [T,\infty)$.  When $\tau \leq T$ Lemmas \ref{lemma:biotsavart} and \ref{lemma:helmholtz} provide the bound $|\omega(\tau)|_\infty\leq C(T)|\tilde{w}(\tau)|_2$, which with (\ref{globalexistencebound0}) implies $(1+|\omega(\tau)|_\infty)\leq C(T)\|\tilde{w}_0\|_m$ for all $\tau \leq T$.  In addition, (\ref{globalexistencebound0}) implies $|\tilde{w}(\tau)|_2\leq C(\|\tilde{w}_0\|_m)$ when $\tau\leq T$.  This concludes the proof.
\end{proof}

\begin{Definition}
Denote by $\Phi(\tilde{w}_0)(\tau)$ the global solution given by Theorem \ref{theorem:VCHEscaledvorticityexistence} corresponding to initial data $\tilde{w}_0\in L^2(m)$.
\end{Definition}

We spend the remaining portion of this section recalling facts about the operator on the spaces $L^2(m)$ useful to our discussion.  The operator $\mathcal{L}$ is studied closely in the appendix of \cite{MR1912106} and the reader is refered there for proofs of the statements which are not immediate.  

The spectrum of $\mathcal{L}$ in the space $L^2(m)$ is
\begin{equation}
\sigma(\mathcal{L})=\{\lambda \in \mathbb{C}  | Re(\lambda)\leq \frac{1-m}{2}\} \cup \{ -\frac{k}{2} | k\in \mathbb{N} \} \notag
\end{equation}

The operator $\mathcal{L}$ forms a strongly continuous semigroup $e^{\tau\mathcal{L}}$ on the space $L^2(m)$.  Set $k= m-2$, then the spectrum of $e^{\mathcal{L}}$ acting on $L^2(m)$ has $k+1$ isolated eigenvalues $\lambda_j = e^{\frac{-j}{2}}$ for $0\leq j\leq k$.  

\begin{Definition}\label{def:subspaces}
Let $X_1 \subset L^2(m)$ be the finite subspace spanned by the eigenvectors associated with the eigenvalues $\lambda_j$, $0\leq j\leq k$ and $X_2 = L^2(m)-X_1$.  Define $P_1$ as the spectral projection onto $X_1$ and $P_2$ the projection onto $X_2$.  
\end{Definition}
\begin{Remark}
Note that $P^k$ and $P^k$ both commute with $e^{\mathcal{L}}$ and are guaranteed to exist because $X_1$ and $X_2$ are closed subspaces of the Hilbert space $L^2(m)$.  This notation supresses the dependence on $m$ but in all applications we shall fix $m$ beforehand so this is well defined.
\end{Remark}

\begin{Lemma}\label{lemma:H3}
The projections $P_1$ and $P_2$ defined above satisfy the following bounds for $\tilde{w}\in L^2(m)$
\begin{align}
\|(H_{\alpha,\tau}e^{\mathcal{L}}P_1)^{-j}w\|_m &\leq C_1 e^{\frac{jk}{2}}\|w\|_m\notag\\
\|(H_{\alpha,\tau}e^{\mathcal{L}}P_2)^{j}w\|_m &\leq C_2 e^{\frac{-j(k+1)}{2}}\|w\|_m\notag
\end{align}
\end{Lemma}
\begin{proof}
This follows from the definitions of the projections.
\end{proof}
This lemma is stated in general but we will work only in the cases $m=2,3$.

Of particular interest to us are the first two eigenvalues $\lambda_0 =1$, $\lambda_1=\frac{-1}{2}$ and the associated eigenvectors $G(\xi)=\frac{1}{4\pi}e^{\frac{-|\xi|}{4}}$ and $F_i= \partial_iG(\xi)=-\frac{\xi}{2}G(\xi)$. The eigenvalue $G$ is particularly important in studying the Navier-Stokes system in 2D, it is a stationary solution of the scaled PDE.  This fact is easily checked by finding the associated velocity field with the Biot-Savart kernel and computing $v^G\cdot\nabla G = 0$.  For reference,
\begin{equation}
v^G(\xi)= \frac{1}{2\pi}\frac{1-e^{\frac{-|\xi|^2}{4}}}{|\xi|^2}\left(\stackrel{-\xi_2}{\xi_1}\right)\notag
\end{equation}

In the scaled VCHE a similar statement is true when we consider $G$ as the filtered and scaled vorticity, this suggests that in describing the behavior of solutions we will also need to consider the following ``un-filtered'' eigenvectors.

\begin{Definition}
\begin{align}
\Gamma(\xi,\tau)=G(\xi)-\alpha^2e^{-\tau}\triangle G(\xi)\notag\\
\Lambda_i(\xi,\tau)=F_i(\xi)-\alpha^2e^{-\tau}\triangle F_i(\xi)\notag
\end{align}
\end{Definition}

Through straight forward calculations one can check 
\begin{align}
v^G\cdot\nabla \Gamma &= 0&\partial_t \Gamma=\mathcal{L}\Gamma &= \alpha^2e^{-\tau}\triangle \Gamma\notag
\end{align}
It is then clear that $\Gamma$ is a solution of the VCHE.  This solution is stationary from the perspective of the filtered flow $\omega$ and plays a similar role as $G$ in the Navier-Stokes equations.

By linearizing the PDE (\ref{PDE:VCHEscaledvorticity}) about the ``fixed'' point $a\Gamma$ ($a\in\mathbb{R}$) we obtain
\begin{align}
\tilde{\psi}_\tau&=\mathcal{L}\tilde{\psi} -a\eta\cdot \nabla \Gamma-a v^G\cdot \nabla \tilde{\psi}&\tilde{\psi}(0)&=\tilde{\psi}_0 \label{PDE:LVCHEscaledvorticity}\\
\eta&=B(\mathcal{H}_{\alpha,\tau}(\tilde{\psi}))\notag
\end{align}
This linear PDE has strong global solutions, a fact which can be established following the steps in Theorem \ref{theorem:VCHEscaledvorticityexistence}; the proof will be omitted here.  Note that $a\Gamma$ is a solution to (\ref{PDE:LVCHEscaledvorticity}).

Letting $v^{F_i}$ denote the velocity field associated with $F_i$, two other important relations are
\begin{align}
v^{F_i}\cdot\nabla \Gamma+v^G\cdot\nabla \Lambda_i &=0&\partial_t \Lambda_i -\mathcal{L}\Lambda_i + \frac{1}{2}\Lambda_i &= 0\notag
\end{align}
The first is quickly checked by differentiating the relation $v^G\cdot\nabla \Gamma = 0$ and the second can be checked directly.  Together, these show $a\Gamma+b_1e^{\frac{-\tau}{2}}\Lambda_1+b_2e^{\frac{-\tau}{2}}\Lambda_2$ is a solution to the linearized equation (\ref{PDE:LVCHEscaledvorticity}).  

Subtracting (\ref{PDE:LVCHEscaledvorticity}) from (\ref{PDE:VCHEscaledvorticity}) we find
\begin{align}
(\tilde{w}-\tilde{\psi})_\tau = \mathcal{L}(\tilde{w}&-\tilde{\psi}) -\omega\cdot\nabla(\tilde{w}-\tilde{\psi})\label{PDE:differencerealation}\\
&-(\omega-av^G)\cdot\nabla\tilde{\psi}+a\eta\cdot\nabla\Gamma\notag
\end{align}
This system will be studied in the following chapters as the foundation for our asymptotic results.

\section{First Order Asymptotic Behavior}
In this section we work in the space $L^2(2)$ where the operator $\mathcal{L}$ has a single isolated eigenvalue $0$ and corresponding eigenvector $G$.  This eigenvector spans the subspace $X_1$ as in Definition \ref{def:subspaces}.  The projection onto this subspace is defined as follows: let $a=\int\tilde{w}_0\,d\xi$, then $P_1(\tilde{w}_0)=aG$.  

Before we proceed it is important to remark that $\int\tilde{w}\,d\xi$ is a conserved quantity under the flow of (\ref{PDE:LVCHEscaledvorticity}) so $P_1(\tilde{w}(\tau))=aG$.  Indeed, since $\mathcal{L}\tilde{w} = \nabla\cdot(\nabla\tilde{w}+\frac{\xi}{2})$,
\begin{equation}
\frac{1}{2}\frac{d}{dt}\int \tilde{w}\,d\xi = \int \nabla\cdot(\nabla\tilde{w}+\frac{\xi}{2}\cdot\tilde{w}-\omega\tilde{w})\,d\xi =0
\end{equation}

Fix $\tilde{w}_0$ and consider the system (\ref{PDE:differencerealation}) with initial conditions $\tilde{w}_0$ and $\tilde{\psi}_0=a\Gamma(\cdot,0)$ where $a=\int\tilde{w}_0\,d\xi$, later we will require $\|\tilde{w}_0\|_m$ to be sufficiently small. With these initial conditions $\tilde{\psi}=a\Gamma$ is a ``stationary'' solution of the linear system (\ref{PDE:LVCHEscaledvorticity}).  Writing $\tilde{f}=\tilde{w}-a\Gamma$, $\phi=B\mathcal{H}_{\alpha,\tau}(\tilde{f})$, (\ref{PDE:differencerealation}) becomes
\begin{equation}
\tilde{f}_\tau = \mathcal{L}\tilde{f}-\omega\cdot\nabla \tilde{f}-a\phi\cdot\nabla\Gamma\label{PDE:VCHEmL}
\end{equation}
The associated integral equation is
\begin{equation}
\tilde{f}=e^{\tau\mathcal{L}}P_2(\tilde{w}_0) - \int_0^\tau e^{\frac{-1}{2}(\tau-\sigma)}\nabla\cdot e^{(\tau-\sigma)\mathcal{L}}(\omega \tilde{f}+a\phi\Gamma)(\sigma)\,d\sigma\notag
\end{equation}

Although it is not obvious at a quick glance, the above system changes with $\tau$.  The dependence on $\tau$ is burried in the $\phi$ term (in the filter relation) and destroys the communitive property of the generated semiflow.  Although communitivity is a very useful property when dealing with semiflows we are able to proceed in a fasion consistent with the Lyapunov-Perreon approach to constructing invariant manifolds, namely to develop a descrete system which represents the solution for integer times and solve this system in a space with suitable decay.  To make this work we need a sequence of semiflows with uniform bounds which can be put together to reconstruct the original flow.  The following system takes us in that direction and we now change our focus to finding properties of $\tilde{f}_n$ which solves the the following system, defined for all $n\in\mathbb{N}$,
\begin{align}
\tilde{f}_n&=e^{\tau\mathcal{L}}\tilde{f}_{n,0} - \int_0^\tau e^{\frac{-1}{2}(\tau-\sigma)}\nabla\cdot e^{(\tau-\sigma)\mathcal{L}}(\omega(n+\sigma) \tilde{f}_n(\sigma)+a\phi_m(\sigma)\Gamma(n+\sigma))\,d\sigma\label{PDE:scrambledVCHE1}\\
\phi_n(\sigma)&=B\mathcal{H}_{\alpha,\sigma+n}(\tilde{f}_n)(\sigma)\ \ \ \ \ \ \ \ \ \tilde{f}(0)=\tilde{f}_0\in X_2\notag
\end{align}
Note that this system has a global solution, again proved analogously to Theorem \ref{theorem:VCHEscaledvorticityexistence} and not included here.  It should be noted here that the above system depends on $\tilde{w}_0$ which we consider fixed when we write down the system ($a$ depends on $\tilde{w}_0$).  Most solutions of this equation will have no meaning but by choosing $\tilde{f}_{0,0}=P_2(\tilde{w}_0)$, then $\tilde{f}_{n,0}=\tilde{f}_{n-1}(1)$ we are able to recover information about the semiflow corresponding to $\tilde{w}_0$.

\begin{Definition}
Let $\Theta_n(\tilde{f}_{n,0})(\tau)$ denote the global solution to the system (\ref{PDE:scrambledVCHE1}).  Sometimes we will write $\Theta$ for $\Theta_0$ when convenient.
\end{Definition}

The parameter $n$ allows us to keep track of our progress in time.  For example, if $0<\sigma\leq 1$ is such that $\tau = n+\sigma$ then
\begin{equation}
\Theta(\tilde{f}_0)(\tau) = \Theta_n(\Theta_{n-1}(\Theta_{n-2}(\cdots)(1))(1))(\sigma)\notag
\end{equation}

That finishes the construction of the semiflows, now we will prove uniform (hold for all $n$) properties.
\begin{Lemma}\label{lemma:H1}
The semiflows $\Theta_n(\tilde{f}_0)(\tau)$, $n\in \mathbb{N}$ are all $C^1$ in $L^2(2)\times \mathbb{R}^+$.  There exists a constant $r_0>0$ (possibly small) and $D>0$ such that for all $\|\tilde{w}_0\|_2<r_0$ and $n\in \mathbb{N}$ the flow $\Theta_n$ satisfies the following Lipschitz property,
\begin{align}\label{bound:H1}
\sup_{0\leq \tau <1} Lip(\Theta_n(\cdot)(\tau)) =D < \infty
\end{align}
This bound holds as $r_0\rightarrow 0$.
\end{Lemma}
\begin{proof}
Since the bound obtained in Lemma \ref{lemma:Helmholtzscaledbound} is independent of $\tau$ we can prove this lemma for all $n\in \mathbb{N}$ at once.  That the semiflow is $C^1$ is a classical result.  Consider the semiflows $\tilde{f}(\tau)=\Theta_n(\tilde{f}_0)(\tau)$ and $\tilde{g}(\tau)=\Theta_n(\tilde{g}_0)(\tau)$ corresponding to initial data $\tilde{f}_0$ and $\tilde{g}_0$ respectively, and the corresponding filtered flows $\phi$ and $\gamma$.  Subtracting we have
\begin{align}
\|\Theta_n(\tilde{f})(\tau)-\Theta_n(\tilde{g})(\tau)\|_2 \leq \|e^{\tau\mathcal{L}}(\tilde{f}_0-\tilde{g}_0)\|_m + I(\tau) + J(\tau)\notag
\end{align}
where
\begin{align}
I(\tau)&= \|\int_0^\tau e^{-\frac{1}{2}(\tau-\sigma)}\nabla \cdot e^{(\tau-\sigma)\mathcal{L}}(\omega(n+\sigma)(\tilde{f}(\sigma)-\tilde{g}(\sigma)))\,d\sigma\|_m\notag\\
J(\tau)&= \|\int_0^\tau e^{-\frac{1}{2}(\tau-\sigma)}\nabla \cdot e^{(\tau-\sigma)\mathcal{L}}((\phi(\sigma)-\gamma(s))a\Gamma(n+\sigma))\,ds\|_m\notag
\end{align}
Similar to the steps in Lemma \ref{lemma:mbilinear} we obtain
\begin{align}
I(\tau)+J(\tau)\leq C(\tau) &\left(\sup_{0\leq \sigma < \tau}\|\tilde{f}(\sigma)-\tilde{g}_2(\sigma)\|_m\right)\cdot\notag\\
&\ \ \ \ \ \cdot\left( \sup_{0\leq \sigma < \tau}(\|\tilde{w}(n+\sigma)\|_m+a\|\Gamma(n+\sigma)\|_m)\right)\notag
\end{align}
with $C(\tau)$ a continuous function and hence bounded on $\tau \in [0,1]$ by a constant $C$.  The bound (\ref{theorem:VCHEscaledvorticitybound})  combined with the definition of $a$ ($a\leq \|\tilde{w}_0\|_2$) allows us to pick $r_0>0$ small enough so that for all $\sigma >0$,
\begin{equation}
\|\tilde{w}(n+\sigma)\|_m+a\|\Gamma(n+\sigma)\|_m < \frac{1}{2C}\notag
\end{equation}
After taking the supremum over $\tau\in[0,1]$ we have
\begin{align}
\sup_{0\leq \tau < 1}\|\Theta_n(\tilde{f})(\tau)-\Theta_n(\tilde{g})(\tau)\|_m 
&\leq \sup_{0\leq \tau < 1}\|e^{\tau\mathcal{L}}(\tilde{f}_0-\tilde{g}_0)\|_m\notag\\
&\leq D\|\tilde{f}_0-\tilde{g}_0\|_m\notag
\end{align}
This is the bound (\ref{bound:H1}).
\end{proof}

\begin{Lemma}\label{lemma:H2}
There exists a constant $r_0>0$ (possibly small) such that for all $\|\tilde{w}_0\|_2<r_0$ and $n\in \mathbb{N}$ the flow $\Theta_n$ can be decomposed as $\Theta_n(\tilde{f}_0)(1) = e^{\mathcal{L}}\tilde{f}_0+R_n(\tilde{f}_0)$ where $R_n(\cdot)$ is Lipschitz as a function from $L^2(2)$ to itself. The Lipschitz constant $Lip(R):=\sup_{n\in\mathbb{N}} Lip(R_n(\cdot))$ can be made arbitrarily small by chosing $r_0$ small and satisfies the following conditions:

(i) For any $\mu\in (0,1/2)$, $r_0$ may be chosen so that for all $i\in\mathbb{N}$,
\begin{align}\label{bound:H42}
\frac{C_1}{1-e^{-\mu}}+\frac{C_2}{e^{-\mu}-e^{-1/2}} <\frac{1}{Lip(R)}
\end{align}

(ii) This bounds holds as $r_0\rightarrow 0$.
\end{Lemma}
\begin{proof}
As in the previous proof, consider the semiflows $\tilde{f}(\tau)=\Theta_n(\tilde{f}_0)(\tau)$ and $\tilde{g}(\tau)=\Theta_n(\tilde{g}_0)(\tau)$.  Define 
\begin{align}
R_n(\tilde{f}_0)&=\Theta_n(\tilde{f}_0)(1) - e^{\mathcal{L}}\tilde{f}_0 \notag\\
&=- \int_0^1 e^{-1/2(1-\sigma)}\nabla\cdot e^{(1-\sigma)\mathcal{L}}(\omega(n+\sigma) \tilde{f}(\sigma)+a\phi(\sigma)\Gamma(n+\sigma))\,d\sigma\notag
\end{align}
As in Lemma \ref{lemma:mbilinear} or the Lemma \ref{lemma:H1}, after adding and subtracting cross terms,
\begin{align}
\|R_n(\tilde{f}_0) - R_n(\tilde{g}_0)\|_m \leq C&\left(\sup_{0\leq \sigma < \tau}\|\tilde{f}(\sigma)-\tilde{g}(\sigma)\|_m\right)\cdot \notag\\
&\ \ \ \ \ \cdot\left(\sup_{0\leq \sigma < \tau}(\|\tilde{w}(n+\sigma)\|_m+a\|\Gamma(n+\sigma)\|_m)\right)\notag
\end{align}
Appealing to (\ref{bound:H1}),
\begin{align}
\|R_n(\tilde{f}_0) - R_n(\tilde{g}_0)\|_m \leq CD&\left(\sup_{0\leq \sigma < \tau}\|\tilde{f}_0-\tilde{g}_0\|_m\right)\cdot \notag\\
&\ \ \ \ \ \cdot\left(\sup_{0\leq \sigma < \tau}(\|\tilde{w}(n+\sigma)\|_m+a\|\Gamma(n+\sigma)\|_m)\right)\notag
\end{align}
This shows 
\begin{equation}
Lip(R_i(\cdot)) \leq CD \left(\sup_{0\leq \sigma < \infty}(\|\tilde{w}(\sigma)\|_m+a\|\Gamma(\sigma)\|_m)\right)\notag
\end{equation}
and, with the help of (\ref{theorem:VCHEscaledvorticitybound}), satisfies the Lipschitz condition.  The bound (\ref{bound:H42}) is established by taking $r_0$ sufficiently small.
\end{proof}

From this point on we assume that $r_0$ is small enough so the conclusions of Lemma \ref{lemma:H1} and \ref{lemma:H2} hold and construct the descrete system.

\begin{Definition}
We call a sequence $\{\tilde{f}_n\}_{n\in\mathbb{N}}\subset L^2(2)$ a positive semiorbit of $\Theta_n$ if it satisfies the relation $\tilde{f}_n = e^{\mathcal{L}}\tilde{f}_{n-1} + R_n(\tilde{f}_{n-1})$.
\end{Definition}

For a given initial position $x\in L^2(m)$, using this definition recursively we find the positive semiorbit $\{\Theta(x)(n)\}_{n\in\mathbb{N}}$ it generates:
\begin{equation}
\Theta(x)(n) = e^{n\mathcal{L}}x +\sum_{k=0}^{n-1} e^{(n-k-1)\mathcal{L}}R_n(\Theta(x)(k))\label{recursiverelation:possemiorbit1}
\end{equation}

We now prove a lemma based on Lemma 3.3 from \cite{MR1472350} which will provide a link between the semiflow $\Theta$ and discrete Lyapunov-Perron equations.
\begin{Lemma}\label{lemma:LPEsemiflow}
With $m=2$, let $X_i$, $P_i$, and $L_i$, for $i=1,2$ be as in Definition \ref{def:subspaces}, pick $\mu\in(0,\frac{1}{2})$ and choose $r_0>0$ to satisfy the conclusions of Lemmas \ref{lemma:H1} and \ref{lemma:H2}.  Let $\{\tilde{f}_n\}_{n\in \mathbb{N}}\subset L^2(2)$ satisfy
\begin{equation}\label{bound:LPEsemiflow}
\limsup_{n\rightarrow\infty}\frac{1}{n}\ln\|\tilde{f}_n\|_2<-\mu
\end{equation}
Then the sequence $\{\tilde{f}_n\}_{n\in \mathbb{N}}$ is a positive semiorbit of $\Theta$ if and only if it satisfies, for all $n\in \mathbb{N}$,
\begin{align}
\tilde{f}_n = e^{n\mathcal{L}}P_2\tilde{f}_0 - \sum_{n\leq j}e^{(n-j-1)\mathcal{L}}P_1R_j(\tilde{f}_j) +\sum_{0\leq j <n}e^{(n-j-1)\mathcal{L}}P_2R_j(\tilde{f}_j)\label{discretePDE:scrambledVCHE1}
\end{align}
\end{Lemma}
\begin{proof}
Using the iterative relation (\ref{recursiverelation:possemiorbit1}) with the projection $P_2$ shows
\begin{equation}
P_2\tilde{f}_n = e^{n\mathcal{L}}P_2\tilde{f}_0 +\sum_{0\leq j <n}e^{(n-j-1)\mathcal{L}}P_2R_j(\tilde{f}_n)\notag
\end{equation}
Likewise,
\begin{equation}
P_1\tilde{f}_n = e^{-m\mathcal{L}}P_1\tilde{f}_0 -\sum_{n\leq j <m}e^{(n-j-1)\mathcal{L}}P_1R_j(\tilde{f}_n)\notag
\end{equation}
We will first show that the right hand side converges as $m\rightarrow \infty$, then it remains only to add these projections together to finish the proof.  Lemma \ref{lemma:H3}  and the Lipschitz property of $R$ allow
\begin{align}
\|\sum_{n\leq j <m}e^{(n-j-1)\mathcal{L}}P^0_1R_j(\tilde{f}_n)\|_m
\leq C_1Lip(R)\sum_{n\leq j <m}\|\tilde{f}_n\|_m\notag
\end{align}
The bound (\ref{bound:LPEsemiflow}) gives convergence of this sum as $m\rightarrow \infty$ as well as the following limit
\begin{align}
\lim_{m\rightarrow\infty}\|e^{-m\mathcal{L}}P^0_1\tilde{f}_n\|_m \leq C_1\|\tilde{f}_n\|_m = 0\notag
\end{align}
\end{proof}

The next step is to show that for any initial data there exists a solution to the system (\ref{discretePDE:scrambledVCHE1}) in a weighted space.  Systems such as (\ref{discretePDE:scrambledVCHE1}) are called Lyapunov-Perron equations and used in the Lyapunov-Perron approach to invariant manifolds.  Theorem 2.1 in \cite{MR1472350} is directly applicable to this system after we introduce the correct weighted space $E^\mu$.  We would like to remark here that the system (\ref{discretePDE:scrambledVCHE1}) does not ``see'' the component of the initial data in $X_1$ but instead picks out the correct component to form a solution with the fast decay rate.  Thus, by solving the system one creates a map $h:X_2\rightarrow X_1$ defined $h(P_2\tilde{f}_0)=P_1\tilde{f}_0$.  This map is important when constructing invariant manifolds and foliations through the Lyapunov-Perron approach (see for example \cite{MR1472350}).  In our construction this map is the zero map.

\begin{Definition}\label{defn:spaceEmu}
Given $\mu\in \mathbb{R}^+$ let $E^\mu_n$ be the Banach space equal to $L^2(2)$ as a vector space but equipped with the norm $|\cdot|_{E^\mu_n}=e^{\mu n}\|\cdot\|_2$.  $E^\mu$ is the sequence space $\tilde{f}=\{\tilde{f}_n\}_{n\in\mathbb{N}}$, $\tilde{f}_n\in E^\mu_n$ equipped with the norm $|\tilde{f}|_{E^\mu}=\sup_{n\in\mathbb{N}}|\tilde{f}_n|_{E^\mu_n}$.
\end{Definition}

\begin{Theorem}\label{theorem:LPEsemiflowexistence}
With $m=2$, let $X_i$, $P_i$, and $L_i$, for $i=1,2$ be as in Definition \ref{def:subspaces}.  Pick $\mu\in(0,\frac{1}{2})$ and choose $r_0>0$ and $\tilde{w}_0$ to satisfy the conclusions of Lemmas \ref{lemma:H1} and \ref{lemma:H2}.  Given initial data $P_2\tilde{f}_0\in X_2$ there exists a unique solution to (\ref{discretePDE:scrambledVCHE1}) in $E^\mu$ with $P_1\tilde{f}_n=0$.
\end{Theorem}
\begin{proof}
Our proof will consist of checking the hypothesis of Theorem 2.1 in \cite{MR1472350} with the space $E^\mu$.  This theorem, which is proved using a fixed point argument, implies the existence and uniqueness statements in our result.

To check (S.1) we note that $\{e^{n\mathcal{L}}P_1\}_{n\in\mathbb{N}}\tilde{f}_0$ is an element of $E^\mu$.  Indeed, from Lemma \ref{lemma:H3}, $e^{\mu n}\|e^{n\mathcal{L}}P_1 \tilde{f}_0\|_2 \leq C_2e^{(\mu-\frac{1}{2}) n}\|\tilde{f}_0\|_2$

To meet hypothesis (R.1) we need to compute the Lipschitz constant associated with the map $R_j$ in the weighted space.  Lemma \ref{lemma:H2} gives the Lipschitz property of $R_j$ as a mapping from $L^2(2)$ to itself with Lipschitz constant $Lip(R_j)<Lip(R)$, this implies $R_j$ is Lipschitz as a mapping $R_j:E^\mu_j\rightarrow E^\mu_n$ with constant $e^{\mu(j-n)}Lip(R)$.  The convergence of the following sums relies on this Lipschitz constant and Lemma \ref{lemma:H3}.
\begin{align}
|S_1|_{E_n^\mu}&=\sum_{n\leq j}e^{\mu n}\|e^{(n-j-1)\mathcal{L}}P_1R_j(\tilde{f}_j)\|_2\notag\\
&\leq C_1\sum_{n\leq j}e^{\mu n}\|R_j(\tilde{f}_j)\|_2\notag\\
&\leq C_1Lip(R)\sum_{n\leq j}e^{\mu n}\|\tilde{f}_j\|_2\notag\\
&\leq C_1Lip(R)\sum_{n\leq j}e^{\mu (n-j)}|\tilde{f}_j|_{E_j^\mu}\notag\\
&\leq C_1Lip(R)\left(\sum_{j\in\mathbb{N}}e^{-\mu j}\right)\left(\sup_{j\in\mathbb{N}}|\tilde{f}_j|_{E_j^\mu}\right)\notag\\
&\leq C_1Lip(R)\left(\frac{1}{1-e^{-\mu}}\right)\left(\sup_{j\in\mathbb{N}}|\tilde{f}_j|_{E_j^\mu}\right)\notag
\end{align}
Similarly,
\begin{align}
|S_2|_{E_n^\mu}&=\sum_{0\leq j <n}\|e^{(n-j-1)\mathcal{L}}P_2R_n(\tilde{f}_j)\|_2\notag\\
&\leq C_2Lip(R)\frac{\sup_{j\in\mathbb{N}}|\tilde{f}_j|_{E_j^\mu}}{e^{-\mu}-e^{-1/2}}\notag
\end{align}
Together with the bound (\ref{bound:H42}) we obtain $|S_1+S_2|_{E_n^\mu}<\sup_{j\in\mathbb{N}}(|\tilde{f}_j|_{E_j^\mu})$ to meet the criteria of (R.1) and grant existence and uniqueness.

The statement $P_1\tilde{f}_n=0$ can be inferred from the ``conservation of mass'' property of (\ref{PDE:VCHEmL}) and the decay implied by $E^\mu$ in the following way.  Using (\ref{PDE:VCHEmL}) we see that $\int \tilde{f}_n\,d\xi$ is a conserved property under the flow $\Theta_n$:
\begin{equation}
\frac{1}{2}\frac{d}{dt}\int \tilde{f}\,d\xi= \int \nabla\cdot(\nabla \tilde{f} + \xi\tilde{f} -\omega\tilde{f}-a\phi\Gamma)\,d\xi =0\notag
\end{equation}
The orthogonal relation $X_1\perp X_2$ allows $\|P_1\tilde{f}_n\|_2 \leq \|\tilde{f}_n\|_2$, and, as $\|P_1\tilde{f}_2\|_m$ is constant for all $n\in \mathbb{N}$, $\{\tilde{f}_n\}\in E^\mu$ implies $P_1\tilde{f}_n = 0$. 
\end{proof}

We are now in a position to prove an existence theorem which will be the basis for the first order decay estimates.

\begin{Theorem}\label{theorem:scrambledexistence1}
With $m=2$, let $X_i$, $P_i$, and $L_i$, for $i=1,2$ be as in Definition \ref{def:subspaces}.  Pick $\mu\in(0,\frac{1}{2})$ and choose $r_0>0$ to satisfy the conclusions of Lemmas \ref{lemma:H1} and \ref{lemma:H2}.  Given $\tilde{w}_0$ such that $\|\tilde{w}_0\|_2\leq r_0$ and initial data $\tilde{f}_0\in X_2$, there exists a unique global solution $\tilde{f}(\tau)\in C^0([0,\infty),L^2(2))$ of (\ref{PDE:scrambledVCHE1}) which satisfies $P_1\tilde{f}(\tau)=0$ and
\begin{equation}\label{bound:scrableddecay1}
\limsup_{\tau\rightarrow\infty}\frac{1}{n}\ln\|\tilde{f}(\tau)\|_2<-\mu
\end{equation}
\end{Theorem}
\begin{proof}
As mentioned earlier, existence of a unique global solution can be argued as in Theorem \ref{theorem:VCHEscaledvorticityexistence}.  As in the previous theorem, conservation of mass and this limit implies $P_1\tilde{f}(\tau)$.   Combining uniqueness, Theorem \ref{theorem:LPEsemiflowexistence} and Lemma \ref{lemma:LPEsemiflow} we can deduce $\limsup_{n\rightarrow\infty}\frac{1}{n}\ln\|\tilde{f}(n)\|_m<-\mu$.  We now apply the Lipschitz property of the semiflow $\Theta_n$ given by Lemma \ref{lemma:H1}.  If $n\in \mathbb{N}$ and $0< \sigma \leq 1$ are such that $\tau = n+\sigma$, then
\begin{align}
\frac{1}{\tau}\ln\|\tilde{f}(\tau)\|_2\leq\frac{1}{\tau}|\ln(D\|\tilde{f}(n)\|_2)\leq\frac{1}{|n|}\ln\|\tilde{f}(n)\|_2 + \frac{1}{\tau}\ln D\notag
\end{align} 
Taking the limit superior as $\tau\rightarrow \infty$ in the above expression finishes the proof.
\end{proof}

We now prove the theorem which was the goal of this section.

\begin{Theorem}Pick $\mu\in(0,\frac{1}{2})$ and choose $r_0>0$ to satisfy the conclusions of Lemmas \ref{lemma:H1} and \ref{lemma:H2}.  Given initial data $\tilde{w}_0$ such that $\|\tilde{w}_0\|_2\leq r_0$, the solution $\tilde{w}(\tau)$ of the scaled VCHE given by Theorem  \ref{theorem:VCHEscaledvorticityexistence} is subject to the following decay estimate:
\begin{equation}
\|\tilde{w}(\tau) - a\Gamma(\tau)\|_2 \leq Ce^{-\mu\tau}\notag
\end{equation}
where $a=\int\tilde{w}\,d\xi$.
\end{Theorem}
\begin{proof}
In the previous theorem take $\tilde{f}_0=P_2\tilde{w}_0$, then $\tilde{f}(\tau)= \tilde{w}(\tau)-a\Gamma(\tau)$.  The decay then follows from (\ref{bound:scrableddecay1}).
\end{proof}

Finally, as a corollary to this theorem we prove Theorem \ref{theorem:decay1}.

\begin{Corollary}\label{cor:gooddecay1}
For any $\mu\in (0,\frac{1}{2})$, there exists a $r_0>0$ so that for any initial data $\tilde{v}_0\in L^2(2)$ such that $\|\tilde{v}_0\|_2\leq r_0$, the solution of (\ref{PDE:VCHEvorticity}) given by Theorem \ref{PDE:VCHEvorticityexistence} satisfies
\begin{equation}
|\tilde{v}(\cdot,t)-a(\Omega(\cdot,t)-\alpha^2\triangle\Omega(\cdot,t))|_p\leq C(1+t)^{-1-\mu+\frac{1}{p}}\notag
\end{equation}
where $a = \int_{\mathbb{R}^2} \tilde{v}_0\,dx$ and
\begin{equation}
\Omega(x,t)= \frac{1}{4\pi(1+t)}e^{\frac{-|x|^2}{4(1+t)}}\notag
\end{equation}
\end{Corollary}
\begin{proof}
This is the result of the previous theorem in unscaled coordinates.  Let
\begin{equation}
\Omega(x,t)= \frac{1}{(1+t)}G(\frac{x}{\sqrt{1+t}})= \frac{1}{4\pi(1+t)}e^{\frac{-|x|^2}{4(1+t)}}\label{function:Omega}
\end{equation}
then,
\begin{equation}
\frac{1}{(1+t)}\Gamma(\frac{x}{\sqrt{1+t}},\ln(1+t)) = \Omega(x,t)-\alpha^2\triangle\Omega(x,t)\notag
\end{equation}
Thanks to the above theorem and the inclusion $L^2(2)\hookrightarrow L^p$ for when $1\leq p \leq 2$,
\begin{align}
|\tilde{v}(\cdot,t)-a(\Omega(\cdot,t)-&\alpha^2\triangle\Omega(\cdot,t))|_p\leq\notag\\
&\ \ \ \ \ \leq(1+t)^{-1+\frac{1}{p}}|\tilde{w}(\cdot,\ln(1+t))-a\Gamma(\cdot,\ln(1+t))|_p\notag\\
&\ \ \ \ \ \leq C(1+t)^{-1+\frac{1}{p}}\|\tilde{w}(\cdot,\ln(1+t))-a\Gamma(\cdot,\ln(1+t))\|_2\notag\\
&\ \ \ \ \ \leq C(1+t)^{-1-\mu+\frac{1}{p}}\notag
\end{align}
\end{proof}

\section{Second Order Asymptotic}

This section is similar in spirit to the previous section and many of the proofs are omitted because they are nearly the same as proofs given in the previous section.  We work in the space $L^2(3)$ where the operator $\mathcal{L}$ has two isolated eigenvalues, $0$ and $\frac{-1}{2}$, and three corresponding eigenvectors, $G$ and $F_i$, $i=1,2$. Together these eigenvalues span the subspace $X_1$ given by Definition \ref{def:subspaces}.  Any $\tilde{w}_0\in L^2(3)$ can be written as $\tilde{w}_0 = aG+b_1F_1+b_2F_2+g$ where $a=\int \tilde{w}_0\,d\xi$, $b_i= \int \xi_i\tilde{w}_0\,d\xi$, and $g\in X_2$.

In addition to the ``conservation of mass'' property $\frac{d}{dt}\int\tilde{w}\,d\xi=0$ discussed in the previous section, solutions to the scaled VCHE (\ref{PDE:VCHEscaledvorticity}) also satisfy the scaled form of conservation of the first moments, $\frac{d}{dt}\int\xi_i\tilde{w}\,d\xi=-\frac{1}{2}\int\xi_i\tilde{w}\,d\xi$. Indeed, let $i=1,2$ and $j\neq i$, then
\begin{align}
\xi_i\mathcal{L}\tilde{w}+\frac{1}{2}\xi_i\tilde{w}=\partial_j(\xi_i\partial_i\tilde{w}+\frac{1}{2}\xi_i^2\tilde{w}-\tilde{w}) +\partial_j(\xi_i\partial_j\tilde{w}+\frac{1}{2}\xi_i\xi_j\tilde{w})\notag
\end{align}
and it is clear that $\int_{\mathbb{R}^2}\xi_i\mathcal{L}\tilde{w}\,d\xi = -\frac{1}{2}\int_{\mathbb{R}^2}\xi_i\tilde{w}\,d\xi$.  For the  non-linear term, note that if $\nabla\cdot\omega=0$ and $\tilde{\omega}=\partial_1 \omega_2-\partial_2\omega_1$ then
\begin{align}
\xi_i\omega\cdot\nabla\tilde{\omega} =\partial_i(\xi_i\omega_i\tilde{\omega}-\omega_i\omega_j) +\partial_j(\xi_i\omega_j\tilde{\omega}+\frac{1}{2}(\omega_i^2-\omega_j^2))\notag
\end{align}
Similarly,
\begin{align}
\xi_i\omega\cdot\nabla\partial_i^2\tilde{\omega} = \partial_i(\xi_i\omega\cdot\nabla\partial_j\tilde{\omega})&- \omega\cdot\nabla\partial_i\tilde{\omega} -\xi_i\partial_i\omega\cdot\nabla\partial_i\tilde{\omega} \notag\\
= \partial_i(\xi_i\omega\cdot\nabla\partial_j\tilde{\omega}) &- \nabla(\omega\partial_i\tilde{\omega}) -\partial_i(\xi_i\partial_i\omega_i\partial_i\tilde{\omega}-\partial_i\omega_i\partial_i\omega_j)\notag\\ 
-\partial_j(\xi_i&\partial_i\omega_j\partial_i\tilde{\omega}+\frac{1}{2}((\partial_i\omega_i)^2-(\partial_i\omega_j)^2))\notag
\end{align}
and
\begin{align}
\xi_i\omega\cdot\nabla\partial_j^2\tilde{\omega} = \partial_j(\xi_i\omega\cdot\nabla\partial_j\tilde{\omega}) &-\partial_i(\xi_i\partial_j\omega_i\partial_j\tilde{\omega}-\partial_j\omega_i\partial_j\omega_j)\notag\\ 
-\partial_j(\xi_i&\partial_j\omega_j\partial_j\tilde{\omega}+\frac{1}{2}((\partial_j\omega_i)^2-(\partial_j\omega_j)^2))\notag
\end{align}
so that $\int_{\mathbb{R}^2}\xi_i\omega\cdot\nabla\tilde{w}\,d\xi = 0$.

Fix $\tilde{w}_0\in L^2(m)$ and consider the system (\ref{PDE:differencerealation}) with initial conditions $\tilde{w}_0$ and $\tilde{\psi}_0=P_1\tilde{w}_0$.  Let $a=\int \tilde{w}_0\,d\xi$, $b_i=\int \xi_i\tilde{w}_0\,d\xi$, then $\tilde{\psi}(\tau)=a\Gamma+e^{\frac{-1}{2}}(b_1\Lambda_1+b_2\Lambda_2)$ is the solution to the linear equation (\ref{PDE:LVCHEscaledvorticity}).  Write $\tilde{f}=\tilde{w}-a\Gamma-e^{\frac{-1}{2}}(b_1\Lambda_1+b_2\Lambda_2)$, $\phi=B\mathcal{H}_{\alpha,\tau}(\tilde{f})$, and after simplification,
\begin{equation}
\tilde{f}_\tau = \mathcal{L}\tilde{f}-\omega\cdot\nabla \tilde{f} -\phi\cdot\nabla\psi-e^{-\tau}(b_1^2v^{F_1}\cdot\nabla\Lambda_1-b_2^2v^{F_2}\cdot\nabla\Lambda_2)\label{PDE:VCHEmL2}
\end{equation}
So, $\tilde{f}$ satisfies the integral equation
\begin{align}
\tilde{f}=e^{\tau\mathcal{L}}&\tilde{f}_0 - \int_0^\tau e^{-1/2(\tau-\sigma)}\nabla\cdot \notag\\
&\cdot e^{(\tau-\sigma)\mathcal{L}}(\omega \tilde{f}+\phi\cdot\nabla\psi+e^{-\tau}(b_1^2v^{F_1}\cdot\nabla\Lambda_1-b_2^2v^{F_2}\cdot\nabla\Lambda_2))(\sigma)\,d\sigma\notag
\end{align}
The key difference between this system and the corresponding system in the previous section is the ``forcing term'' $e^{-\tau}(b_1^2v^{F_1}\cdot\nabla\Lambda_1-b_2^2v^{F_2}\cdot\nabla\Lambda_2)$.  In general such terms may add difficulty but this one decays at a rate fast enough to be little more then a notational nuisance.

As in the previous section, fix an initial condition $\tilde{w}_0$ and for any $n\in\mathbb{N}$ consider the system,
\begin{align}
\tilde{f}=e^{\tau\mathcal{L}}\tilde{f}_0 - \int_0^\tau e^{-1/2(\tau-\sigma)}&\nabla\cdot e^{(\tau-\sigma)\mathcal{L}}(\omega(n+\sigma) \tilde{f}(\sigma)+a\phi(\sigma)\psi(n+\sigma))\,d\sigma\label{PDE:scrambledVCHE2}\\
- \int_0^\tau e^{-1/2(\tau-\sigma)}\nabla\cdot &e^{(\tau-\sigma)\mathcal{L}}(e^{-(n+\sigma)}(b_1^2v^{F_1}\cdot\nabla\Lambda_1-b_2^2v^{F_2}\cdot\nabla\Lambda_2)(n+\sigma))\,d\sigma\notag\\
\phi(\sigma)&=B\mathcal{H}_{\alpha,\sigma+\tau}\notag(\tilde{f})(\sigma) \ \ \ \ \ \ \ \ \ \tilde{f}(0)=\tilde{f}_0\in X_2\notag
\end{align}

\begin{Definition}
Let $\Psi_n(\tilde{f}_0)(\tau)$ denote the global solution to the system (\ref{PDE:scrambledVCHE2}).
\end{Definition}

\begin{Lemma}\label{lemma:H12}
For each $n\in \mathbb{N}$, the semiflows $\Psi_n(\tilde{f}_0)(\tau)$ is $C^1$ in $L^2(2)\times \mathbb{R}^+$.  There exists a constant $r_0>0$ (possibly small) and $D>0$ such that for all $\|\tilde{w}_0\|_3<r_0$ and $n\in \mathbb{N}$ the flow $\Psi_n$ satisfies the following Lipschitz property,
\begin{align}\label{bound:H12}
\sup_{0\leq \tau <1} Lip(\Psi_n(\cdot)(\tau)) =D < \infty
\end{align}
This bound holds as $r_0\rightarrow 0$.
\end{Lemma}
\begin{proof}
The proof is nearly identical to the proof of Lemma \ref{lemma:H1}
\end{proof}

\begin{Lemma}\label{lemma:H22}
There exists a constant $r_0>0$ (possibly small) such that for all $\|\tilde{w}_0\|_3<r_0$ and $n\in \mathbb{N}$ the flow $\Psi_n$ can be decomposed as $\Psi_n(\tilde{f}_0)(1) = e^{\mathcal{L}}\tilde{f}_0+R_n(\tilde{f}_0)+S_n$ where $R_n(\cdot)$ is Lipschitz as a function from $L^2(3)$ to itself and $S_n\in L^2(3)$. The Lipschitz constant $Lip(R):=\sup_{n\in\mathbb{N}} Lip(R_n(\cdot))$ can be made arbitrarily small and satisfies the following conditions:

(i) For any $\mu\in (1/2,1)$, $r_0$ may be chosen so that for all $n\in\mathbb{N}$,
\begin{align}\label{bound:H422}
\frac{C_1}{e^{-1/2}-e^{-\mu}}+\frac{C_2}{e^{-\mu}-e^{-1}} <\frac{1}{Lip(R)}
\end{align}

(ii) This bounds all hold as $r_0\rightarrow 0$.
\end{Lemma}
\begin{proof}
Define
\begin{equation}
S_n=- \int_0^\tau e^{-1/2(\tau-\sigma)}\nabla\cdot e^{(\tau-\sigma)\mathcal{L}}(e^{-(n+\sigma)}(b_1^2v^{F_1}\cdot\nabla\Lambda_1-b_2^2v^{F_2}\cdot\nabla\Lambda_2)(n+\sigma))\,d\sigma\notag
\end{equation}
and $R_n(\tilde{f}_0)=\Psi_n(\tilde{f}_0)(1) - e^{\mathcal{L}}\tilde{f}_0 - S_n$.  To see that $S_n\in L^2(3)$ is similar to the proof of Lemma \ref{lemma:mbilinear}:
\begin{align}
\|S_n\|_3\leq e^{-n}(b_1\|v^{F_1}\|_3\|\Lambda_1\|_3+b_2\|v^{F_2}\|_3\|\Lambda_2\|_3)\label{bound:Sn}
\end{align}
The remaining statements in this proof follow as in the proof of Lemma \ref{lemma:H2}.
\end{proof}

\begin{Lemma}\label{lemma:LPEsemiflow2}
With $m=3$, let $X_i$, $P_i$, and $L_i$, for $i=1,2$ be as in Definition \ref{def:subspaces}.  Pick $\mu\in(\frac{1}{2},1)$ and choose $r_0>0$ to satisfy the conclusions of Lemmas \ref{lemma:H12} and \ref{lemma:H22}.  Let $\{\tilde{f}_n\}_{n\in \mathbb{N}}\subset L^2(3)$ satisfy
\begin{equation}\label{bound:LPEsemiflow2}
\limsup_{n\rightarrow\infty}\frac{1}{n}\ln\|\tilde{f}_n\|_3<-\mu
\end{equation}
Then the sequence $\{\tilde{f}_n\}_{n\in \mathbb{N}}$ is a positive semiorbit of $\Psi$ if and only if it satisfies, for all $n\in \mathbb{N}$,
\begin{align}
\tilde{f}_n = e^{n\mathcal{L}}P_2\tilde{f}_0 &- \sum_{n\leq j}e^{(n-j-1)\mathcal{L}}P_1(R_j(\tilde{f}_j)+S_j)+\sum_{0\leq j <n}e^{(n-j-1)\mathcal{L}}P_2(R_j(\tilde{f}_j)+S_j) \label{discretePDE:scrambledVCHE2}
\end{align}
\end{Lemma}
\begin{proof}
Proof is nearly identical to Lemma \ref{lemma:LPEsemiflow}
\end{proof}

For the purpose of this section, we take $E^\mu$ as in Definition \ref{defn:spaceEmu} but with the base space of $L^2(3)$ instead of $L^2(2)$.

\begin{Theorem}\label{theorem:LPEsemiflowexistence2}
With $m=3$, let $X_i$, $P_i$, and $L_i$, for $i=1,2$ be as in Definition \ref{def:subspaces}.  Pick $\mu\in(\frac{1}{2},1)$ and choose $r_0>0$ to satisfy the conclusions of Lemmas \ref{lemma:H12} and \ref{lemma:H22}.  Given initial data $P_2\tilde{f}_0\in X_2$ there exists a unique solution to (\ref{discretePDE:scrambledVCHE2}) in $E^\mu$ with $P_1\tilde{f}_n=0$.
\end{Theorem}
\begin{proof}
Again, the proof is nearly identical to Theorem \ref{theorem:LPEsemiflowexistence}, the main difference is the ``forcing term'' $S_n$ which we include with the linear term and results in one extra check to satisfy hypothesis of Theorem 2.1 in \cite{MR1472350}: show that $\{S_n\}_{n\in\mathbb{N}}\in E^\mu$.  This follows from the bound (\ref{bound:Sn}) and $\mu<1$: $e^{\mu n}\|S_n\|_3\leq C e^{-(1-\mu) n}\leq C$
\end{proof}

\begin{Theorem}\label{theorem:scrambledexistence2}
With $m=3$, let $X_i$, $P_i$, and $L_i$, for $i=1,2$ be as in Definition \ref{def:subspaces}.  Pick $\mu\in(\frac{1}{2},1)$ and choose $r_0>0$ to satisfy the conclusions of Lemmas \ref{lemma:H12} and \ref{lemma:H22}.  Given $\tilde{w}_0$ such that $\|\tilde{w}_0\|_3\leq r_0$ and initial data $\tilde{f}_0\in X_2$, there exists a unique global solution $\tilde{f}(\tau)\in C^0([0,\infty),L^2(2))$ of (\ref{PDE:scrambledVCHE1}) which satisfies $P_1\tilde{f}(\tau)$ and
\begin{equation}
\limsup_{\tau\rightarrow\infty}\frac{1}{n}\ln\|\tilde{f}(\tau)\|_3<-\mu
\end{equation}
\end{Theorem}
\begin{proof}
Proof is nearly identical to Theorem \ref{theorem:scrambledexistence1}.
\end{proof}

\begin{Theorem}With $m=3$, let $X_i$, $P_i$, and $L_i$, for $i=1,2$ be as in Definition \ref{def:subspaces}.  Pick $\mu\in(\frac{1}{2},1)$ and choose $r_0>0$ to satisfy the conclusions of Lemmas \ref{lemma:H1} and \ref{lemma:H2}.  Given initial data $\tilde{w}_0$ such that $\|\tilde{w}_0\|_3\leq r_0$, the solution $\tilde{w}(\tau)$ of the scaled VCHE given by Theorem  \ref{theorem:VCHEscaledvorticityexistence} is subject to the following decay estimate:
\begin{equation}
\|\tilde{w}(\tau) - a\Gamma(\tau)-e^{\frac{-\tau}{2}}(b_1\Lambda_1+b_2\Lambda_2)\|_2 \leq Ce^{-\mu\tau}\notag
\end{equation}
where $a=\int\tilde{w}\,d\xi$ and $b_i=\int\xi_i\tilde{w}\,d\xi$.
\end{Theorem}
\begin{proof}
In the previous theorem take $\tilde{f}_0=P_2\tilde{w}_0$, then $\tilde{f}(\tau)= \tilde{w}(\tau)-a\Gamma(\tau)-e^{\frac{-\tau}{2}}(b_1\Lambda_1+b_2\Lambda_2)$.  The decay then follows from (\ref{bound:scrableddecay1}).
\end{proof}

The following corollary proves Theorem \ref{theorem:decay2}

\begin{Corollary}
For any $\mu\in (\frac{1}{2},1)$, there exists a $r_0>0$ so that for any initial data $\tilde{v}_0\in L^2(3)$ with $\|\tilde{v}_0\|_3\leq r_0$ the solution of (\ref{PDE:VCHEvorticity}) given by Theorem \ref{PDE:VCHEvorticityexistence} satisfies
\begin{equation}
|\tilde{v}(\cdot,t)-a(\Omega(\cdot,t)-\alpha^2\triangle\Omega(\cdot,t))-\sum_{i=1,2}b_i(\partial_i\Omega(x,t)-\alpha^2\triangle\partial_i\Omega(x,t))|_p\leq C(1+t)^{-1-\mu+\frac{1}{p}}\notag
\end{equation}
where $a=\int\tilde{w}\,d\xi$ and $b_i=\int\xi_i\tilde{w}\,d\xi$ and $\Omega$ is defined by (\ref{function:Omega}).
\end{Corollary}
\begin{proof}
This is the result of the previous theorem in unscaled coordinates.  Let $\Omega(x,t)= \frac{1}{(1+t)}G(\frac{x}{\sqrt{1+t}})$, then,
\begin{align}
\frac{1}{(1+t)}\Gamma\left(\frac{x}{\sqrt{1+t}},\ln(1+t)\right) = \Omega(x,t)-\alpha^2\triangle\Omega(x,t)\notag\\
e^{-\frac{\tau}{2}}\frac{1}{(1+t)}\Lambda_i\left(\frac{x}{\sqrt{1+t}},\ln(1+t)\right) = \partial_i\Omega(x,t)-\alpha^2\triangle\partial_i\Omega(x,t)\notag
\end{align}
The rest follows as in Corollary \ref{cor:gooddecay1}.
\end{proof}

\bibliographystyle{plain}
\bibliography{vche-2dasym}

\begin{thebibliography}{10}

\bibitem{MR1308857}
M.~Ben-Artzi.
\newblock Global solutions of two-dimensional {N}avier-{S}tokes and {E}uler
  equations.
\newblock {\em Arch. Rational Mech. Anal.}, 128(4):329--358, 1994.

\bibitem{BjorlandSchonbek06}
C.~Bjorland and M.~E. Schonbek.
\newblock On questions of decay and existence for the viscous camassa-holm
  equations.
\newblock {\em Submitted}, 2006.

\bibitem{MR1308858}
H.~Brezis.
\newblock Remarks on the preceding paper by {M}. {B}en-{A}rtzi: ``{G}lobal
  solutions of two-dimensional {N}avier-{S}tokes and {E}uler equations''
  [{A}rch.\ {R}ational {M}ech.\ {A}nal.\ {\bf 128} (1994), no.\ 4, 329--358;
  {MR}1308857 (96h:35148)].
\newblock {\em Arch. Rational Mech. Anal.}, 128(4):359--360, 1994.

\bibitem{MR1381974}
E.~A. Carlen and M.~Loss.
\newblock Optimal smoothing and decay estimates for viscously damped
  conservation laws, with applications to the {$2$}-{D} {N}avier-{S}tokes
  equation.
\newblock {\em Duke Math. J.}, 81(1):135--157 (1996), 1995.
\newblock A celebration of John F. Nash, Jr.

\bibitem{MR1274542}
A.~Carpio.
\newblock Asymptotic behavior for the vorticity equations in dimensions two and
  three.
\newblock {\em Comm. Partial Differential Equations}, 19(5-6):827--872, 1994.

\bibitem{MR1472350}
X.~Chen, J.~K. Hale, and B.~Tan.
\newblock Invariant foliations for {$C\sp 1$} semigroups in {B}anach spaces.
\newblock {\em J. Differential Equations}, 139(2):283--318, 1997.

\bibitem{MR1837927}
C.~Foias, D.~D. Holm, and E.~S. Titi.
\newblock The {N}avier-{S}tokes-alpha model of fluid turbulence.
\newblock {\em Phys. D}, 152/153:505--519, 2001.
\newblock Advances in nonlinear mathematics and science.

\bibitem{MR1878243}
C.~Foias, D.~D. Holm, and E.~S. Titi.
\newblock The three dimensional viscous {C}amassa-{H}olm equations, and their
  relation to the {N}avier-{S}tokes equations and turbulence theory.
\newblock {\em J. Dynam. Differential Equations}, 14(1):1--35, 2002.

\bibitem{MR1912106}
T.~Gallay and C.~E. Wayne.
\newblock Invariant manifolds and the long-time asymptotics of the
  {N}avier-{S}tokes and vorticity equations on {$\mathbb{R}^2$}.
\newblock {\em Arch. Ration. Mech. Anal.}, 163(3):209--258, 2002.

\bibitem{MR953819}
Y.~Giga and T.~Kambe.
\newblock Large time behavior of the vorticity of two-dimensional viscous flow
  and its application to vortex formation.
\newblock {\em Comm. Math. Phys.}, 117(4):549--568, 1988.

\bibitem{MR1627802}
D.~D. Holm, J.~E. Marsden, and T.~S. Ratiu.
\newblock The {E}uler-{P}oincar\'e equations and semidirect products with
  applications to continuum theories.
\newblock {\em Adv. Math.}, 137(1):1--81, 1998.

\bibitem{htlans}
D.~D. Holm and E.~S. Titi.
\newblock Computational models of turbulence: The lans-$\alpha$ model and the
  role of global analysis.
\newblock {\em SIAM News}, 38, 2005.

\bibitem{MR2031580}
A.~A. Ilyin and E.~S. Titi.
\newblock Attractors for the two-dimensional {N}avier-{S}tokes-{$\alpha$}
  model: an {$\alpha$}-dependence study.
\newblock {\em J. Dynam. Differential Equations}, 15(4):751--778, 2003.

\bibitem{MR1270113}
T.~Kato.
\newblock The {N}avier-{S}tokes equation for an incompressible fluid in {${\bf
  R}\sp 2$} with a measure as the initial vorticity.
\newblock {\em Differential Integral Equations}, 7(3-4):949--966, 1994.

\bibitem{MR1853633}
J.~E. Marsden and S.~Shkoller.
\newblock Global well-posedness for the {L}agrangian averaged {N}avier-{S}tokes
  ({LANS}-{$\alpha$}) equations on bounded domains.
\newblock {\em R. Soc. Lond. Philos. Trans. Ser. A Math. Phys. Eng. Sci.},
  359(1784):1449--1468, 2001.
\newblock Topological methods in the physical sciences (London, 2000).

\bibitem{MR916761}
H.~Osada.
\newblock Diffusion processes with generators of generalized divergence form.
\newblock {\em J. Math. Kyoto Univ.}, 27(4):597--619, 1987.

\end{thebibliography}
\end{document}